%% file: SIAM.tex
\documentclass[a4paper,american,reqno]{amsart}

\def\MyAMSFlag{1}
\def\MyInformsFlag{0}

\input{setup}

\bibliography{references}

\title[Decision Problems in Multilevel Linear Programming]{Decision Problems in \\ Multilevel Linear Programming}
\author{Nagisa Sugishita}
\author{Margarida Carvalho}
\address[Nagisa Sugishita, Margarida Carvalho]{André-Aisenstadt Pavillon, 2920 Tour Road, Montreal, Quebec H3T 1N8, Canada}
\email{nagisa.sugishita@hec.ca}

\date{\today}

\begin{document}

\begin{abstract}
\input{abstract}
\end{abstract}

\maketitle

\input{body.tex}

\printbibliography

\clearpage

\appendix

\input{appendix.tex}

\end{document}

%% file: setup.tex
\usepackage{babel}
\usepackage[T1]{fontenc}
\usepackage[utf8]{inputenc}
\usepackage{siunitx}
\usepackage[dvipsnames,table]{xcolor}
\usepackage{amsmath}
\usepackage{amsfonts}
\usepackage{amssymb}
\usepackage{mathtools}
\usepackage{booktabs}
\usepackage{multirow}
\usepackage{enumitem}
\usepackage{mathdots}
\setlist[enumerate,1]{label=(\roman*)}

\if\MyAMSFlag1
\usepackage[style=numeric,
            giveninits=true,
            maxbibnames=100,
            maxcitenames=2,
            isbn=false,
            doi=true,
            url=false,
            isbn=false,
            backend=biber]{biblatex}
\newcommand{\citet}[1]{\citeauthor{#1}~\cite{#1}}

\newenvironment{MyProof}[1][Proof.]{\begin{proof}[#1]}{\end{proof}}
\newcommand{\MyQED}{}
\fi
\usepackage[colorlinks,
            citecolor=blue,
            urlcolor=blue,
            linkcolor=blue]{hyperref} % should always be the last package

\if\MyAMSFlag1
\newtheorem{theorem}{Theorem}

\newtheorem{lemma}[theorem]{Lemma}

\fi

\if\MyInformsFlag1

\newenvironment{MyProof}[1][Proof]{\begin{proof}{#1.}}{\end{proof}}

\newcommand{\MyQED}{\hfill\halmos}
\fi

\input{macros}

%% file: macros.tex
\newcommand{\MacroColor}[1]{#1}

\newcommand{\NonNegativeIntegers}{\mathbb{Z}_{\ge 0}}
\newcommand{\PositiveIntegers}{\mathbb{Z}_{> 0}}

\newcommand{\ZeroVector}{\mathbf{0}}
\newcommand{\OneVector}{\mathbf{1}}

\DeclareMathOperator*{\Argmin}{\MacroColor{argmin}}

\newcommand{\MaximumEntrySize}{\MacroColor{\sigma}}
\newcommand{\OptValFunc}{\MacroColor{v}}
\newcommand{\OptSolSet}{\MacroColor{\mathcal{S}}}
\newcommand{\FeasSolSet}{\MacroColor{\mathcal{F}}}

\newcommand{\ComplexityClassWrapper}[1]{\MacroColor{\ensuremath{#1}}}
\newcommand{\ComplexityClassText}[1]{\textup{#1}}

\newcommand{\ClassNP}{\ComplexityClassWrapper{\ComplexityClassText{NP}}}
\newcommand{\ClassCoNP}{\ComplexityClassWrapper{\ComplexityClassText{coNP}}}

\newcommand{\ClassPiP}[1]{\ComplexityClassWrapper{\Pi^{p}_{#1}}}
\newcommand{\ClassSigmaP}[1]{\ComplexityClassWrapper{\Sigma^{p}_{#1}}}
\newcommand{\ClassDeltaP}[1]{\ComplexityClassWrapper{\Delta^{p}_{#1}}}

\newcommand{\ClassPSPACE}{\ComplexityClassWrapper{\ComplexityClassText{PSPACE}}}

\newcommand{\ComplexityProblemWrapper}[1]{\MacroColor{\ensuremath{#1}}}
\newcommand{\ComplexityProblemText}[1]{\textup{#1}}

\newcommand{\DecisionProblemQThreeSAT}[1]{\ComplexityProblemWrapper{#1\ComplexityProblemText{-QSAT}}}

\newcommand{\EqTagDecisionProblemOnKLP}[2]{\ComplexityProblemWrapper{#1\ComplexityProblemText{-#2}}}

\newcommand{\DecisionProblemModifierComp}[1]{\ensuremath{\overline{#1}}}

\newcommand{\DecisionProblemFeas}[1]{\hyperref[LabelDecisionProblemFeas]{\EqTagDecisionProblemOnKLP{#1}{FEAS}}}
\newcommand{\DecisionProblemVal}[1]{\hyperref[LabelDecisionProblemVal]{\EqTagDecisionProblemOnKLP{#1}{VAL}}}
\newcommand{\DecisionProblemUnb}[1]{\hyperref[LabelDecisionProblemUnb]{\EqTagDecisionProblemOnKLP{#1}{UNB}}}
\newcommand{\DecisionProblemRecognize}[1]{\hyperref[LabelDecisionProblemRecognize]{\EqTagDecisionProblemOnKLP{#1}{RECOG}}}
\newcommand{\DecisionProblemAttain}[1]{\hyperref[LabelDecisionProblemAttain]{\EqTagDecisionProblemOnKLP{#1}{ATAIN}}}

\newcommand{\NVariables}{\MacroColor{n}}
\newcommand{\NConstraints}{\MacroColor{m}}

\newcommand{\SATBooleanFormula}{\MacroColor{F}}
\newcommand{\SATQBF}{\MacroColor{H}}
\newcommand{\SATQBFSet}[1]{\MacroColor{\mathcal{B}_{#1}}}

\newcommand{\SATBoolVariable}{\MacroColor{\nu}}

\newcommand{\EqTagGeneralKLP}[3]{\MacroColor{\ensuremath{\mathrm{LP}_{\ensuremath{#1}}^{\ensuremath{#2}}(\ensuremath{#3})}}}
\newcommand{\EqrefGeneralKLP}[3]{\textup{(}\text{\hyperref[LabelGeneralKLP]{$\EqTagGeneralKLP{#1}{#2}{#3}$}\textup{)}}}

\newcommand{\EqTagJeroslowPenaltyBefore}[3]{\MacroColor{\ensuremath{\mathrm{J}_{\ensuremath{#1}}^{\ensuremath{#2}}(\ensuremath{#3})}}}
\newcommand{\EqrefJeroslowPenaltyBefore}[3]{\textup{(}\text{\hyperref[LabelJeroslowPenaltyBefore]{$\EqTagJeroslowPenaltyBefore{#1}{#2}{#3}$}\textup{)}}}

\newcommand{\EqTagJeroslowPenaltyAfter}[3]{\MacroColor{\ensuremath{\overline{\mathrm{J}}_{\ensuremath{#1}}^{\ensuremath{#2}}(\ensuremath{#3})}}}
\newcommand{\EqrefJeroslowPenaltyAfter}[3]{\textup{(}\text{\hyperref[LabelJeroslowPenaltyAfter]{$\EqTagJeroslowPenaltyAfter{#1}{#2}{#3}$}\textup{)}}}

\newcommand{\EqTagQThreeSATReducedToKLPLeadersProblem}[3]{\MacroColor{\ensuremath{\textup{QLP}_{\ensuremath{#1}}^{\ensuremath{#2}}\ensuremath{#3}}}}
\newcommand{\EqTagQThreeSATReducedToKLPFollowersProblem}[3]{\MacroColor{\ensuremath{\textup{QLP}_{\ensuremath{#1}}^{\ensuremath{#2}}\ensuremath{#3}}}}

\newcommand{\EqRefQThreeSATReducedToKLPLeadersProblem}[3]{\textup{(}\text{\hyperref[eq:LabelQThreeSATReducedToKLPLeadersProblem]{$\EqTagQThreeSATReducedToKLPLeadersProblem{#1}{#2}{#3}$}\textup{)}}}
\newcommand{\EqRefQThreeSATReducedToKLPFollowersProblem}[3]{\textup{(}\text{\hyperref[eq:LabelQThreeSATReducedToKLPLeadersProblem]{$\EqTagQThreeSATReducedToKLPFollowersProblem{#1}{#2}{#3}$}\textup{)}}}

\newcommand{\HomogeneousLemmaScalar}{\MacroColor{\lambda}}

%% file: abstract.tex
% \begin{abstract}
We study the computational complexity of decision problems in $k$-level linear programming (LP). Seminal work by Jeroslow establishes that determining whether the optimal objective value of a $k$-level LP is at least as good as a given threshold is \ClassSigmaP{k-1}-hard. In this paper, we demonstrate the matching upper bound and thereby prove that this problem is \ClassSigmaP{k-1}-complete. To this end, we show that the feasible region of a $k$-level LP can be expressed as a union of sets defined by weak and strict linear inequalities. Moreover, we show that the decision of the unboundedness is \ClassSigmaP{k-1}-complete. Finally, we discuss the extension of our results to the mixed-binary cases. In short, this work closes lasting open questions in multilevel programming.

%% file: body.tex
% Associate editor
% % Juan Pablo Vielma, Ivana Ljubic, Willem-Jan van Hoeve
% Reviewer
% Shriram, Ted,, Oleg, Meta theorems, Buccheim, Basu

\section{Introduction}\label{sec:introduction}

% context
A multilevel programming problem is an optimization problem where some of the variables are constrained to be in the optimal solution set of another optimization problem.
Since the pioneering papers of~\citet{bracken1973mathematical} and~\citet{candler1977multi}, research on this topic has grown substantially.
In particular, advances in bilevel programming have led to a rich body of literature that studies both its theoretical foundations and practical solution methods.
See \cite{carvalho2025integer,dempe2020bilevel} for recent developments.
There is also emerging interest in trilevel programming due to its flexible modeling capability, such as the critical node problem~\cite{nabli2022complexity}, resource allocation problems~\cite{cassidy1971efficient}, fortification-interdiction
games~\cite{tomasaz2024completeness}, and cyber security~\cite{liu2015trilevel}.
More applications of multilevel programming can be found in~\cite{migdalas2013multilevel,vicente1994bilevel}.

% what we study in the paper
Among the models of multilevel programming explored in the literature, we focus on the so-called optimistic variant.
This means that if there are multiple optimal solutions for a lower-level problem, the one most favorable to the upper-level problems is chosen. More concretely, in this work, we study the computational complexity of decision problems in multilevel linear programming (LP).
A common approach to analyzing the complexity of optimization problems is to consider the following decision problem: determining whether the optimal objective value of a given instance is less than or equal to a specified threshold, assuming a minimization problem. 
This is often referred to as the decision version of the optimization problem.

% bilevel programming literature
The computational complexity of bilevel LP has received widespread attention in the literature.
In a seminal paper, \citet{Jeroslow1985} showed that the decision version of bilevel LP is \ClassNP{}-hard.
\citet{ben1990computational} gave another proof of \ClassNP{}-hardness.
Later, \citet{hansen1992new} showed that the problem remains \ClassNP{}-hard even when restricted to the min-max case.
Recently, \citet{Buchheim2023} proved the inclusion of this decision problem in \ClassNP{}, implying its \ClassNP{}-completeness.
There are also studies that consider other types of decision problems associated with bilevel LP.
\citet{RodriguesEtAl2024} showed that deciding whether the problem is unbounded (i.e., whether there is a sequence of feasible solutions whose objective values diverge to $-\infty$) is NP-complete.
\citet{VicenteEtAl1994} investigated the recognition (decision) of local optimality of a given point and showed that it is \ClassCoNP{}-hard.
For further discussion about local optimal solutions, see \citet{ProkopyevRalphs2024}.
Recently, Sugishita and Carvalho~\cite{sugishita2025complexitybilevellinearprogramming} showed that the decision version of bilevel LP remains \ClassNP{}-complete even if the input instance is restricted to have a single upper-level variable.

Multilevel LP with more than two levels remains relatively unexplored.
In the aforementioned work, \citet{Jeroslow1985} showed that the decision version of $k$-level LP is \ClassSigmaP{k-1}-hard.
\citet{blair1992computational} presented an alternative argument building on the knapsack problem.
\citet{RodriguesEtAl2024} showed that deciding the unboundedness of $k$-level LP is \ClassSigmaP{k-1}-hard, assuming the upper-level problems have linear constraints (known as coupling or linking constraints).
% In a survey, \citet{deng1998complexity} claimed that $k$-level LP is $F\Delta^{\text{p}}_k$-complete, citing an unpublished reference.
% To the best of our knowledge, no formal proof for this claim exists.
\citet{Buchheim2023} briefly discussed extending his argument for the inclusion of the decision version of bilevel LP in \ClassNP{} to the multilevel setting, suggesting that the decision version of $k$-level LP is in \ClassSigmaP{k-1} (Remark 2 in the reference).
However, we believe that extending his argument to general multilevel LPs is nontrivial due to subtle issues that arise for $k$-level LPs with $k \ge 3$, as discussed in Online Supplement~A.
\citet{benson1989structure} and \citet{bard1985geometric} studied geometric properties of multilevel LP.

In this paper, we strengthen the results of \citet{Jeroslow1985} and \citet{blair1992computational} by establishing matching upper bounds, namely that the decision version of $k$-level LP belongs to $\ClassSigmaP{k-1}$. 
In particular, we prove $\ClassSigmaP{k-1}$-completeness. 
Our approach is based on a novel structural analysis of $k$-level LP.
We show that the feasible region of a rational $k$-level LP instance can be expressed as a union of sets defined by rational weak and strict inequalities of polynomial (encoding) size. 
This extends the result of \citet{BasuEtAl2021}, which shows that the feasible set of a bilevel LP instance is a union of polyhedra.
In particular, our result implies that every rational feasible $k$-level LP instance admits a rational feasible solution of polynomial size.
Moreover, we establish $\ClassSigmaP{k-1}$-completeness of the problem of deciding unboundedness for $k$-level LP. 
A related hardness result was recently obtained in \cite{RodriguesEtAl2024}. 
However, in contrast to their reduction, our construction works under weaker assumptions, namely, without upper-level constraints. 
Finally, we briefly discuss how our arguments can be extended to prove, for the first time, $\ClassSigmaP{k}$-completeness of the decision version of $k$-level mixed-binary LP.

We emphasize that our result concerning the computational complexity of the unboundedness of $k$-level LP is in a strong sense. 
Specifically, our proof of $\ClassSigmaP{k-1}$-hardness of the decision problem for the unboundedness of $k$-level LP relies on a reduction that produces a $k$-level LP instance with integer coefficients whose magnitudes are polynomial in the size of the original input. 
The same holds true for the proof of $\ClassSigmaP{k-1}$-hardness for the decision version of $k$-level LP, as established by \citet{Jeroslow1985} (in contrast to the proof by \citet{blair1992computational}, which relies on instances with coefficients of exponential magnitude).
In analogy with the notion of strong $\ClassNP$-hardness, we therefore use the term strong $\ClassSigmaP{k-1}$-hardness to emphasize that the hardness results persist even when all numerical data are polynomially bounded. 
More precisely, this means that, unless $\ClassDeltaP{k-1} = \ClassSigmaP{k-1}$, there exists no Turing machine with oracle access to $\ClassSigmaP{k-2}$ that decides the above decision problems in time polynomial in both the encoding size of the instance and the maximum absolute value of its coefficients (assuming all data is scaled to integers).

\textbf{Paper Structure.}
This paper is organized as follows.
In Section~\ref{sec:MainResults}, we introduce the definition of $k$-level LP and our main results, \ClassSigmaP{k-1}-completeness of both the decision version of $k$-level LP and the decision version of unboundedness of $k$-level LP.
Section~\ref{sec:ProofOfHardness} discusses their \ClassSigmaP{k-1}-hardness, while Section~\ref{sec:ProofOfMembership} shows their membership in \ClassSigmaP{k-1}.
Section~\ref{sec:extensions_and_limitations} provides a brief overview of the extensions and limitations of our argument, with particular emphasis on the computational complexity of $k$-level mixed-binary LP.
Section~\ref{sec:Conclusions} concludes the paper.

\textbf{Conventions.}
We say that a vector, matrix, or optimization instance is \emph{rational} if and only if all its data are rational numbers.
Given an optimization instance $P$, we use $\OptValFunc(P)$, $\FeasSolSet(P)$, and $\OptSolSet(P)$ to denote the optimal objective value, the set of feasible solutions, and the set of optimal solutions, respectively.
We write $\NonNegativeIntegers{} = \{0, 1, \ldots, \}$ and $\PositiveIntegers{} = \{1, 2, \ldots, \}$.
We use $\ZeroVector$ (respectively, $\OneVector$) to denote the all-zeros (respectively, all-ones) vector, with the dimension clear from the context.
We use subscripts to denote elements of a vector: for a vector $v$, the $i$th element is denoted by $v_i$.
To reduce clutter, we sometimes refer to column vectors inline; given $u \in \mathbb{R}^{d_1}$ and $v \in \mathbb{R}^{d_2}$, we use
the notation $(u, v)$ to denote the column vector
$
\begin{pmatrix} u \\ v \end{pmatrix} \in \mathbb{R}^{d_1 + d_2}.
$

\section{Main Results}
\label{sec:MainResults}

Let $k \in \PositiveIntegers{}$, and let $n_l \in \PositiveIntegers{}$, $m_l \in \NonNegativeIntegers{}$ for each $l = 1, \ldots, k$.
For each $l = 1, \ldots, k$ and $i = 1, \ldots, k$, let $A_{l i} \in \mathbb{Q}^{\NConstraints_l \times \NVariables_i}$, $b_l \in \mathbb{Q}^{\NConstraints_l}$ and $c_{l i} \in \mathbb{Q}^{\NVariables_i}$.
% ObjTruncation
We collectively write $A = \{ A_{li} : l = 1, \ldots, k, i = 1, \ldots, k \}$, $b = \{ b_{l} : l = 1, \ldots, k \}$, and $c = \{ c_{l i} : l = 1, \ldots, k, i = l, \ldots, k \}$.
A $k$-level LP instance $(A, b, c)$ consists of a sequence of $k$  optimization problems.
% ObjTruncation
Let the first player choose a decision variable $x_1 \in \mathbb{R}^{n_1}$, followed by the second player choosing $x_2 \in \mathbb{R}^{n_2}$, and so on, until the $k$-th player selects $x_k \in \mathbb{R}^{n_k}$. The objective function of the $l$-th player is given by $\sum_{i = l}^k c_{l i}^{\top} x_i$. 
% For the $l$-th player with $l \ge 2$, the inclusion of the variables $x_{l'}$ for $l' < l$ in the objective does not affect the optimization, since these values are given as fixed constants for the $l$-th player. 
Formally, the first player's problem~\eqref{LabelGeneralKLP} is defined as
\begin{gather}
\inf_{x_1, \ldots, x_k}
\left\{ 
\sum_{i = 1}^k c_{1 i}^{\top} x_i 
: 
\displaystyle
\sum_{i = 1}^k A_{1 i} x_i \ge b_1, 
\begin{pmatrix}
x_2 \\ \vdots \\ x_k
\end{pmatrix}
\in 
\mathcal{S}\eqref{LabelGeneralKLPk_1Level}
\right\},
\label{LabelGeneralKLP}
\tag{$\EqTagGeneralKLP{k}{1}{A, b, c}$}
\end{gather}
where the second player's problem~\eqref{LabelGeneralKLPk_1Level} is
% ObjTruncation
\begin{gather}
\inf_{x_2, \ldots, x_k}
\left\{ \sum_{i = 2}^{k} c_{2 i}^{\top} x_i : 
\displaystyle
\sum_{i = 1}^{k} A_{2 i} x_i \ge b_2, 
\begin{pmatrix}
x_3 \\ \vdots \\ x_k
\end{pmatrix}
\in 
\mathcal{S}(\EqTagGeneralKLP{k}{3}{x_1, x_2, A, b, c})
\right\},
\label{LabelGeneralKLPk_1Level}
\tag{$\EqTagGeneralKLP{k}{2}{x_1, A, b, c}$}
\end{gather}
with the $l$-th player's problem $\EqTagGeneralKLP{k}{l}{x_1, \ldots, x_{l - 1}, A, b, c}$, $3 \le l \le k - 1$, defined analogously, down to the $k$-th player's problem~\eqref{eq:1Level}:
% ObjTruncation
\begin{gather}
\inf_{x_k} 
\left\{ 
c_{k k}^{\top} x_k 
:
\sum_{i = 1}^k 
A_{k i} x_i \ge b_k
\right\}.
\label{eq:1Level}
\tag{$\EqTagGeneralKLP{k}{k}{x_1, \ldots, x_{k - 1}, A, b, c}$}
\end{gather}
We identify the $k$-level LP instance $(A, b, c)$ with the first player's problem~\eqref{LabelGeneralKLP}. Accordingly, the optimal objective value of the $k$-level LP instance $(A, b, c)$ is given by $\OptValFunc\eqref{LabelGeneralKLP}$, and a solution $(x_1, \ldots, x_k)$ is feasible for this instance if and only if $(x_1, \ldots, x_k) \in \FeasSolSet\eqref{LabelGeneralKLP}$.

We consider the following decision problems:
\begin{align}
&
\parbox{0.85\textwidth}{Given a rational $k$-level LP instance $(A, b, c)$ and a rational number $t$, is there a feasible solution whose objective value is less than or equal to $t$?
% is the optimal objective value less than or equal to $t$?
}
\tag{\EqTagDecisionProblemOnKLP{k}{VAL}}
\label{LabelDecisionProblemVal}
\\[0.5em]
&
\parbox{0.85\textwidth}{Given a rational $k$-level LP instance $(A, b, c)$, is it unbounded from below, i.e., $\OptValFunc\eqref{LabelGeneralKLP} = -\infty$?
}
\tag{\EqTagDecisionProblemOnKLP{k}{UNB}}
\label{LabelDecisionProblemUnb}
\end{align}

To capture the computational complexity of \DecisionProblemVal{k} and \DecisionProblemUnb{k} precisely, we give particular consideration to instances satisfying the following restrictions.
\begin{enumerate}[label=(C\arabic*)]
\item
\label{Condition1}
Except for the last player's problem~\eqref{eq:1Level}, there are no linear constraints;
\item 
\label{Condition2}
The linear constraints in the last player's problem~\eqref{eq:1Level} include variable bounds $\ZeroVector \le x_i \le \OneVector$ for all $i = 1, \ldots, k$;
\item 
\label{Condition3}
All entries in $A$, $b$ and $c$ are integers between $-n$ and $n$, where $n = n_1 + \cdots + n_k$.
\end{enumerate}
% As is standard, the specific polynomial bound in~\ref{Condition3}, here $n$, is not essential.
% In fact, once a similar result is established for $n^{\mathcal{O}(1)}$, the same result follows for $n$ by appropriately padding the instances with additional variables.

The main results of this paper are stated below.
\begin{theorem}
\label{theorem:ComplexityOfMultilevelProgramming}
\mbox{} % // TODO
\begin{enumerate}
\item For any $k$, the decision problem \DecisionProblemVal{k} is \ClassSigmaP{k - 1}-complete.
Moreover, this remains true even when the input is restricted to satisfy conditions~\ref{Condition1}--\ref{Condition3}.
\item For any $k$, the decision problem \DecisionProblemUnb{k} is \ClassSigmaP{k - 1}-complete. Moreover, this remains true even when the input is restricted to satisfy conditions~\ref{Condition1} and~\ref{Condition3}.
\end{enumerate}
\end{theorem}

\citet{Jeroslow1985} showed that \DecisionProblemVal{k} is \ClassSigmaP{k - 1}-hard under conditions~\ref{Condition1}--\ref{Condition3}.
In this paper, we show \ClassSigmaP{k - 1}-hardness of \DecisionProblemUnb{k} under conditions~\ref{Condition1} and \ref{Condition3} (Section~\ref{sec:ProofOfHardness}), as well as the membership of \DecisionProblemVal{k} and \DecisionProblemUnb{k} in \ClassSigmaP{k - 1} (Section~\ref{sec:ProofOfMembership}), implying their completeness.

\section{Proof of Hardness}
\label{sec:ProofOfHardness}

In this section, we prove the hardness.

\begin{lemma}
\label{lemma:HardnessOfDecisionProblemValAndDecisionProblemUnb}
\mbox{}
\begin{enumerate}
\item
\label{lemma:HardnessOfDecisionProblemValAndDecisionProblemUnb:one}
The decision problem \DecisionProblemVal{k} is \ClassSigmaP{k - 1}-hard. 
Moreover, this remains true even when the input is restricted to satisfy conditions~\ref{Condition1}--\ref{Condition3}.
\item
\label{lemma:HardnessOfDecisionProblemValAndDecisionProblemUnb:two}
The decision problem \DecisionProblemUnb{k} is \ClassSigmaP{k - 1}-hard.
Moreover, this remains true even when the input is restricted to satisfy conditions~\ref{Condition1} and~\ref{Condition3}.
\end{enumerate}
\end{lemma}

Assertion \ref{lemma:HardnessOfDecisionProblemValAndDecisionProblemUnb:one} in Lemma~\ref{lemma:HardnessOfDecisionProblemValAndDecisionProblemUnb} is shown by \citet{Jeroslow1985} (and by \citet{blair1992computational} without condition~\ref{Condition3}).
The proof is based on the reduction from the quantified satisfiability problem as defined below.
For $k \in \PositiveIntegers{}$, let $\SATQBFSet{k}$ be the set of closed quantified Boolean formulae (QBF) in prenex normal form of the form
$$
\exists \SATBoolVariable_{(1)} \forall \SATBoolVariable_{(2)} \exists \SATBoolVariable_{(3)} \cdots Q \SATBoolVariable_{(k)} \, [\SATBooleanFormula(\SATBoolVariable_{(1)}, \SATBoolVariable_{(2)},  \SATBoolVariable_{(3)}, \ldots, \SATBoolVariable_{(k)}) = 1],
$$
where the set of Boolean variables $\SATBoolVariable_{(l)}$ for $l=1,\ldots,k$ are pairwise disjoint and nonempty, $\SATBooleanFormula$ is a Boolean formula without quantifiers, and $Q$ is $\exists$ if $k$ is odd and $\forall$ otherwise.
Formally, $\SATQBFSet{k}$ is a set of encodings of such QBF.
However, throughout this work, we identify an encoding with the QBF it represents, and we say that $\SATQBF \in \SATQBFSet{k}$ is true (holds) if the encoded QBF hods.
For $k=1,2,\ldots$, we define $\DecisionProblemQThreeSAT{k} = \{\SATQBF \in \SATQBFSet{k} : \text{$\SATQBF$ is true}\}$.
The decision problem \DecisionProblemQThreeSAT{k} is known to be \ClassSigmaP{k}-complete~\cite{wrathall1976complete}.
For each $k \in \PositiveIntegers{}$, \citet{Jeroslow1985} demonstrated a transformation from $\SATQBF \in \SATQBFSet{k}$ into a $(k + 1)$-level LP instance $(A^{\SATQBF}, b^{\SATQBF}, c^{\SATQBF})$ such that 
\begin{equation}
\OptValFunc\EqrefGeneralKLP{k + 1}{1}{A^{\SATQBF}, b^{\SATQBF}, c^{\SATQBF}} = 
\begin{cases}
0 & \text{ if } \SATQBF \in \DecisionProblemQThreeSAT{k}, \\
1 & \text{ otherwise.}
\end{cases}
\label{Eq:OptimalObjectiveValueOfJeroslowsInstance}
\end{equation}
The $(k+1)$-level instance $(A^{\SATQBF}, b^{\SATQBF}, c^{\SATQBF})$ satisfies conditions~\ref{Condition1} and~\ref{Condition2}. 
Moreover, the instance has integer cost coefficients of polynomial magnitude. 
Consequently, by introducing appropriately many dummy variables, one can additionally ensure that condition~\ref{Condition3} is satisfied, which implies assertion~\ref{lemma:HardnessOfDecisionProblemValAndDecisionProblemUnb:one} in Lemma~\ref{lemma:HardnessOfDecisionProblemValAndDecisionProblemUnb}.

To prove assertion~\ref{lemma:HardnessOfDecisionProblemValAndDecisionProblemUnb:two} in Lemma~\ref{lemma:HardnessOfDecisionProblemValInstance}, we use the following results.

\begin{lemma}
\label{lemma:homogeneous-klp-without-linking-constraints}
Let $(A, b, c)$ be a $k$-level LP instance~\eqref{LabelGeneralKLP} satisfying condition~\ref{Condition1}.
Then, for any $\lambda \in \mathbb{R}$ such that $\lambda > 0$,
$\mathcal{F}\EqrefGeneralKLP{k}{1}{A, \lambda b, c} = \lambda \mathcal{F}\eqref{LabelGeneralKLP}$.
\end{lemma}

\begin{MyProof}
We prove the following: for any $l=1,\ldots,k$, $x_1,\ldots,x_{l-1}$, and $\lambda > 0$,
\begin{align}
\HomogeneousLemmaScalar \OptSolSet\EqrefGeneralKLP{k}{l}{x_1, \ldots, x_{l-1}, A, b, c}
= \OptSolSet\EqrefGeneralKLP{k}{l}{\HomogeneousLemmaScalar x_1, \ldots, \HomogeneousLemmaScalar x_{l-1}, A, \HomogeneousLemmaScalar b, c}.
\label{eq:lemma-homogenuity}
\end{align}
We use induction starting from $l=k$.
We only show the induction steps, since the one for the base case is similar.
Suppose \eqref{eq:lemma-homogenuity} holds for some integer $l+1$.
Let $S = \OptSolSet\EqrefGeneralKLP{k}{l}{\HomogeneousLemmaScalar x_1, \ldots, \HomogeneousLemmaScalar x_{l-1}, A, \HomogeneousLemmaScalar b, c}$.
Then,
% ObjTruncation
\begin{align*}
S &= \arg\inf_{x_l, \ldots, x_k} \left\{ \sum_{i=l}^k c_{l i}^{\top} x_i : \begin{pmatrix}x_{l+1}\\\vdots\\x_k\end{pmatrix} \in \OptSolSet\EqrefGeneralKLP{k}{l+1}{\HomogeneousLemmaScalar x_1, \ldots, \HomogeneousLemmaScalar x_{l-1}, x_l, A, \HomogeneousLemmaScalar b, c} \right\}.
\end{align*}
By induction hypothesis,
% ObjTruncation
\begin{align*}
S &= \arg\inf_{x_l, \ldots, x_k} \left\{ \sum_{i=l}^k c_{l i}^{\top} x_i : \begin{pmatrix}x_{l+1}\\\vdots\\x_k\end{pmatrix} \in \HomogeneousLemmaScalar \OptSolSet\EqrefGeneralKLP{k}{l+1}{x_1, \ldots, x_{l-1}, x_l/\HomogeneousLemmaScalar, A, b, c} \right\} \\
&= \arg\inf_{x_l, \ldots, x_k} \left\{ \sum_{i=l}^k c_{l i}^{\top} x_i/\HomogeneousLemmaScalar : \begin{pmatrix}x_{l+1}/\HomogeneousLemmaScalar\\\vdots\\x_k/\HomogeneousLemmaScalar\end{pmatrix} \in \OptSolSet\EqrefGeneralKLP{k}{l+1}{x_1, \ldots, x_{l-1}, x_l/\HomogeneousLemmaScalar, A, b, c} \right\},
\end{align*}
which is equal to $\HomogeneousLemmaScalar \OptSolSet\EqrefGeneralKLP{k}{l}{x_1, \ldots, x_{l-1}, A, b, c}$.
Thus, \eqref{eq:lemma-homogenuity} holds for $l$ as well.
\MyQED
\end{MyProof}

% ----------------------------------

The next lemma shows that constraints not involving the followers’ variables can be ``moved'' to the next player’s problem.

\begin{lemma}
\label{lemma:constraint-forwarding-induction-step}
Let $(A, b, c)$ be a $k$-level LP instance.
Let $p \le k - 1$ and
$$
A_{p i} = \begin{pmatrix}
A_{p i}'
\\
A_{p i}''
\end{pmatrix},
\ \ 
b_{p} = \begin{pmatrix}
b_{p}'
\\
b_{p}''
\end{pmatrix}
$$
for each $i = 1, \ldots, k$, such that $A_{p i}' = 0$ for all $i = p + 1, \ldots, k$.
Define a $k$-LP instance $(\hat{A}, \hat{b}, \hat{c})$ by $\hat{c} = c$, $\hat{A}_{l i} = A_{l i}$, $\hat{b}_{l} = b_{l}$ for all $l = 1, \ldots, k$, $i = 1, \ldots, k$, such that $l \not= p, p + 1$, and
\begin{align*}
\hat{A}_{p i}
=
A_{p i}'', \ \ 
\hat{b}_p = b_p'',
\qquad
\ \ 
\hat{A}_{p + 1 \, i}
=
\begin{pmatrix}
A_{p i}' \\
A_{p + 1 \, i}
\end{pmatrix}
, \ \ 
\hat{b}_{p + 1} = 
\begin{pmatrix}
b_p ' \\ b_{p + 1}
\end{pmatrix},
\quad
i = 1, \ldots, k.
\end{align*}
Then, $\FeasSolSet\EqrefGeneralKLP{k}{1}{A, b, c} = \FeasSolSet\EqrefGeneralKLP{k}{1}{\hat{A}, \hat{b}, \hat{c}}$.
\end{lemma}

\begin{MyProof}
We only show the case with $p \le k - 2$, since the case with $p = k - 1$ is analogous.

Instance $(\hat{A}, \hat{b}, \hat{c})$ is obtained by ``moving'' constraints in the $p$-th player's problem to the $(p + 1)$-th player's problem.
In particular, for $l = p + 2, \ldots, k$, the $l$-th player's problems of $(A, b, c)$ and $(\hat{A}, \hat{b}, \hat{c})$ are identical.
The $(p + 1)$-th player's problems of the two instances are
% ObjTruncation
\begin{align}
\label{eq:prop:constraint-forwarding-p-minus-one-level-1}
% \RFeasSolSet{k}{p - 1}(x_p, \ldots, x_k, \hat{A}, \hat{b}, \hat{c})
\inf_{x_{p + 1}, \ldots, x_k} \left\{
\sum_{i = p + 1}^{k} c_{p + 1 \, i}^T x_i 
:
\begin{array}{l}
\displaystyle
\sum_{i = 1}^{k} A_{p + 1 \, i} x_i \ge b_{p + 1}, \\
(x_{p + 2}, \ldots, x_k)
\in
\OptSolSet\EqrefGeneralKLP{k}{p + 2}{x_1, \ldots, x_{p + 1}, A, b ,c}
\end{array}
\right\}
\end{align}
and
% ObjTruncation
\begin{align}
\label{eq:prop:constraint-forwarding-p-minus-one-level-2}
% \RFeasSolSet{k}{p - 1}(x_p, \ldots, x_k, \hat{A}, \hat{b}, \hat{c})
\inf_{x_{p + 1}, \ldots, x_k} \left\{
\sum_{i = p + 1}^{k} c_{p + 1 \, i}^T x_i 
:
\begin{array}{l}
\displaystyle
\sum_{i = 1}^{p} A_{p i}' x_i \ge b_p', 
\
\displaystyle
\sum_{i = 1}^{k} A_{p + 1 \, i} x_i \ge b_{p + 1}, \\
% \begin{pmatrix}
% x_1 \\ \vdots \\ x_{p - 2}
% \end{pmatrix}
(x_{p + 2}, \ldots, x_k)
\in
\OptSolSet\EqrefGeneralKLP{k}{p + 2}{x_1, \ldots, x_{p + 1}, A, b ,c}
\end{array}
\right\}.
\end{align}
Let
$
Q
=
\{
% \begin{pmatrix}
% x_{p} \\
% \vdots \\
% x_{k}
% \end{pmatrix}
(x_1, \ldots, x_p)
:
\sum_{i = 1}^{p} A_{p i}' x_{i} \ge b_{p}'
\}.
$
For any $(x_1, \ldots, x_p) \in Q$, \eqref{eq:prop:constraint-forwarding-p-minus-one-level-1} and \eqref{eq:prop:constraint-forwarding-p-minus-one-level-2} coincide, while for any $(x_1, \ldots, x_p) \not\in Q$, \eqref{eq:prop:constraint-forwarding-p-minus-one-level-2} is infeasible.
Thus, $\OptSolSet\EqrefGeneralKLP{k}{p + 1}{x_1, \ldots, x_p, \hat{A}, \hat{b}, \hat{c}} = \OptSolSet\EqrefGeneralKLP{k}{p+1}{x_1, \ldots, x_p, A, b, c}$ for all $ (x_1, \ldots, x_p) \in Q$ and $\OptSolSet\EqrefGeneralKLP{k}{p + 1}{x_1, \ldots, x_p, \hat{A}, \hat{b}, \hat{c}} = \emptyset$ otherwise.
% \begin{align*}
% \OptSolSet\EqrefGeneralKLP{k}{p + 1}{x_1, \ldots, x_p, \hat{A}, \hat{b}, \hat{c}}
% &=
% \begin{cases}
% \OptSolSet\EqrefGeneralKLP{k}{p+1}{x_1, \ldots, x_p, A, b, c}, 
% & \quad \forall (x_1, \ldots, x_p) \in Q, \\
% \emptyset, 
% & \quad \forall (x_1, \ldots, x_p) \not\in Q.
% \end{cases}
% \end{align*}
It follows that for any $x_1, \ldots, x_{p - 1}$,
\begin{align*}
&
\FeasSolSet\EqrefGeneralKLP{k}{p}{x_1, \ldots, x_{p - 1}, \hat{A}, \hat{b}, \hat{c}}
\\
&=
\left\{
\begin{pmatrix}
x_p \\ \vdots \\ x_k
\end{pmatrix}
:
\displaystyle
\sum_{i = 1}^{k} A_{p i}'' x_i \ge b_p'', 
\
\begin{pmatrix}
x_1 \\ \vdots \\ x_p
\end{pmatrix}
\in
Q,
\
\begin{pmatrix}
x_{p + 1} \\ \vdots \\ x_k
\end{pmatrix}
\in
\OptSolSet\EqrefGeneralKLP{k}{p+1}{x_1, \ldots, x_p, \hat{A}, \hat{b}, \hat{c}}
\right\}
\\
&=
\left\{
\begin{pmatrix}
x_p \\ \vdots \\ x_k
\end{pmatrix}
:
\displaystyle
\sum_{i = 1}^{k} A_{p i}'' x_i \ge b_p'', 
\
\begin{pmatrix}
x_1 \\ \vdots \\ x_p
\end{pmatrix}
\in
Q,
\
\begin{pmatrix}
x_{p + 1} \\ \vdots \\ x_k
\end{pmatrix}
\in
\OptSolSet\EqrefGeneralKLP{k}{p+1}{x_1, \ldots, x_p, A, b, c}
\right\}
\\
&=
\FeasSolSet\EqrefGeneralKLP{k}{p}{x_1, \ldots, x_{p - 1}, A, b, c}.
\end{align*}
Therefore, for any $x_1, \ldots, x_{p - 1}$, the $p$-th player's problems coincide, and hence $\OptSolSet\EqrefGeneralKLP{k}{p}{x_1, \ldots, x_{p - 1}, \hat{A}, \hat{b}, \hat{c}} = \OptSolSet\EqrefGeneralKLP{k}{p}{x_1, \ldots, x_{p - 1}, A, b, c}$.
Thus, inductively, we conclude that the feasible sets and the optimal solution sets of the $l$-th player's problems are the same for all $l = p - 1, \ldots, 1$.
\MyQED
\end{MyProof}

% ----------------------------------

\begin{lemma}
\label{lemma:constraint-forwarding}
Let $(A, b, c)$ be a rational $k$-level LP instance satisfying
$$
A_{l i} = \begin{pmatrix}
A_{l i}'
\\
A_{l i}''
\end{pmatrix},
\ \ 
b_{l} = \begin{pmatrix}
b_{l}'
\\
b_{l}''
\end{pmatrix},
\quad
\forall l = 1, \ldots, k -1, i = 1, \ldots, k,
$$
such that $A_{l i}' = 0$ for all $l = 1, \ldots, k - 1$, $i = l + 1, \ldots, k$.
Define a $k$-LP instance $(\hat{A}, \hat{b}, \hat{c})$ by $\hat{A}_{l i} = A_{l i}''$, $\hat{b}_{l} = b_{l}''$ for all $l = 1, \ldots, k - 1$, $i = 1, \ldots, k$, and
\begin{equation*}
\begin{pmatrix}
\hat{A}_{k 1}
& 
\hat{A}_{k 2}
& 
\cdots
&
\hat{A}_{k k}
\end{pmatrix}
=
\begin{pmatrix}
A'_{11} & & & \\
\vdots & \ddots & &  \\
A'_{k - 1 \, 1} &  & A'_{k - 1 \, k - 1} & \\
A_{k 1} & \cdots & A_{k k - 1} &  A_{k k} 
\end{pmatrix}, \ \ 
\hat{b} = \begin{pmatrix}
b_{1}' \\
\vdots \\
b_{k - 1}' \\
b_{k} \\
\end{pmatrix},
\ \ 
\hat{c} = c.
\end{equation*}
Then, 
$\FeasSolSet\EqrefGeneralKLP{k}{1}{\hat{A}, \hat{b}, \hat{c}} = \FeasSolSet\EqrefGeneralKLP{k}{1}{A, b, c}$ and $\OptValFunc\EqrefGeneralKLP{k}{1}{\hat{A}, \hat{b}, \hat{c}} = \OptValFunc\EqrefGeneralKLP{k}{1}{A, b, c}$.
\end{lemma}

\begin{MyProof}
It follows from Lemma~\ref{lemma:constraint-forwarding-induction-step}.
\MyQED
\end{MyProof}

Now, we are ready to prove Lemma~\ref{lemma:HardnessOfDecisionProblemValAndDecisionProblemUnb}.

\begin{MyProof}[Proof of Lemma~\ref{lemma:HardnessOfDecisionProblemValAndDecisionProblemUnb}]
Assertion~\ref{lemma:HardnessOfDecisionProblemValAndDecisionProblemUnb:one} is shown in \citet{Jeroslow1985}, so we need only show assertion~\ref{lemma:HardnessOfDecisionProblemValAndDecisionProblemUnb:two}.
The assertion holds for $k=1$ from the polynomial solvability of LP~\cite{Schrijver1998}.
In the following, we assume $k\ge2$.
Let $\SATQBF \in \SATQBFSet{k - 1}$ and let $(\bar{A}^{\SATQBF}, \bar{b}^{\SATQBF}, \bar{c}^{\SATQBF})$ be the $k$-level LP instance satisfying conditions~\ref{Condition1}--\ref{Condition3} such that $\OptValFunc\EqrefGeneralKLP{k}{1}{\bar{A}^{\SATQBF}, \bar{b}^{\SATQBF}, \bar{c}^{\SATQBF}} = -1$ if $\SATQBF \in \DecisionProblemQThreeSAT{(k-1)}$ and $0$ otherwise.
This instance is obtained by subtracting 1 (or, equivalently, by subtracting a variable fixed to 1) from the first player’s objective function in the $k$-level LP instance constructed by \citet{Jeroslow1985} (see \eqref{Eq:OptimalObjectiveValueOfJeroslowsInstance}).
Now, consider the following $k$-level LP instance, which has an extra variable $\lambda_1 \in \mathbb{R}$ in the first player's problem:
\begin{equation}
v'(\SATQBF) = 
\inf_{x_1, \ldots, x_k, \lambda_1 \ge 1}
\left\{ 
\sum_{i = 1}^k (c_{1 i}^{\SATQBF})^{\top} x_i 
: 
\begin{pmatrix} x_2 \\ \vdots \\ x_k \end{pmatrix}
\in 
\mathcal{S}\EqrefGeneralKLP{k}{2}{x_1, \bar{A}^{\SATQBF}, \lambda_1 \bar{b}^{\SATQBF}, \bar{c}^{\SATQBF}}
\right\}.
\label{lemma:HardnessOfDecisionProblemValInstance}
\end{equation}
By Lemma~\ref{lemma:homogeneous-klp-without-linking-constraints}, we have
$$
v'(\SATQBF) 
= 
\inf_{\lambda_1 \ge 1} \lambda_1 \OptValFunc\EqrefGeneralKLP{k}{1}{\bar{A}^{\SATQBF}, \bar{b}^{\SATQBF}, \bar{c}^{\SATQBF}}
=
\begin{cases}
-\infty & \text{ if $\SATQBF \in \DecisionProblemQThreeSAT{(k-1)}$,} \\
0 & \text{ otherwise,}
\end{cases}
$$
implying $\DecisionProblemQThreeSAT{(k-1)} \le_l \DecisionProblemUnb{k}$.
In light of Lemma~\ref{lemma:constraint-forwarding}, one can rewrite \eqref{lemma:HardnessOfDecisionProblemValInstance} as a $k$-level LP instance satisfying conditions~\ref{Condition1} and~\ref{Condition3} (the variable bound $\lambda_1 \ge 1$ can be moved to the last follower's problem).
This completes the proof.
\MyQED
\end{MyProof}

\section{Proof of Membership}
\label{sec:ProofOfMembership}

In this section, we prove that \DecisionProblemVal{k} and \DecisionProblemUnb{k} are in \ClassSigmaP{k - 1}.
Section~\ref{sec:PiecewiseLinearFunctions} studies properties of rational generalized polyhedra and rational piecewise linear functions, which are defined therein.
Section~\ref{sec:FeasibleSetOfMultilevelLP} examines the feasible set of a rational $k$-level LP instance.
By establishing that the value functions of the lower-level problems are rational piecewise linear functions, we show that the feasible set of a $k$-level LP is a union of rational generalized polyhedra.
Section~\ref{sec:ProofOfInclusionInSigmaPKMinusOne:MainAnalysis} then shows that \DecisionProblemVal{k} and \DecisionProblemUnb{k} are in \ClassSigmaP{k - 1}, based on the results from the preceding sections.

As already remarked in Section~\ref{sec:introduction}, Buchheim~\cite{Buchheim2023} established that \DecisionProblemVal{2} belongs to \ClassNP{}.
He further suggested that his argument could be generalized to show that \DecisionProblemVal{k} lies in \ClassSigmaP{k - 1}, and provided a brief sketch of this reasoning (Remark 2 in the reference).
However, we believe that generalizing his argument to multilevel LPs is nontrivial due to subtle issues that arise for $k$-level LPs when $k \ge 3$, as discussed in Online Supplement~A.

\subsection{Generalized Polyhedra and Piecewise Linear Functions}
\label{sec:PiecewiseLinearFunctions}

We say that a set $P \subseteq \mathbb{R}^n$ is a \emph{rational generalized polyhedron} if and only if it can be written as $P = \{ x \in \mathbb{R}^n : A x \ge a, B x > b \}$ for some rational matrices $A$, $B$ and rational vectors $a$, $b$ of conforming dimensions.
In this case, we say that the inequalities $A x \ge a$ and $B x > b$ define $P$.
% We say that $P$ is a rational generalized polyhedron if it is defined by a system of rational linear inequalities.

It is well known that if a rational generalized polyhedron is non-empty, it contains a rational point of size polynomial in the number of variables and the maximum size of the inequalities defining the polyhedron.
% Below, we establish an analogous result for a rational generalized polyhedron.

\begin{lemma}[Schrijver~\cite{Schrijver1998}]
\label{lemma:estimate-of-weak-and-strong-linear-inequality}
Let $n \in \PositiveIntegers{}$ and $m, m' \in \NonNegativeIntegers{}$.
Let $Q = \{x \in \mathbb{R}^{\NVariables} : A x \ge a, B x > b\}$ be a nonempty generalized polyhedron, where $A \in \mathbb{Q}^{\NConstraints \times \NVariables}$, $B \in \mathbb{Q}^{\NConstraints' \times \NVariables}$, $a \in \mathbb{Q}^{\NConstraints}$, $b \in \mathbb{Q}^{\NConstraints'}$.
Let $\MaximumEntrySize$ be the maximum size among the entries in $A$, $a$, $B$ and $b$.
% p134
% where each inequality in Ax 5b has size at most p.
% P44
% each row of the matrix [A b] has size at most p
\begin{enumerate}
\item 
\label{lemma:estimate-of-weak-and-strong-linear-inequality:one}
The closure of $Q$ is given by $\{ x \in \mathbb{R}^{\NVariables} : A x \ge a, B x \ge b \}$.
\item 
\label{lemma:estimate-of-weak-and-strong-linear-inequality:two}
Set $Q$ contains a rational point of size polynomially bounded by $\MaximumEntrySize$ and $n$.
% \item 
% \label{lemma:estimate-of-weak-and-strong-linear-inequality:three}
% If $\inf \{ c^{\top} x : x \in Q \}$ is attainable, it is attained by a rational point of size polynomial in $\MaximumEntrySize$ and $n$.
\end{enumerate}
\end{lemma}

The projection of a rational generalized polyhedron is itself a rational generalized polyhedron. Moreover, the size of the inequalities defining the projection can be bounded polynomially in the number of variables and in the maximum size of the coefficients in the inequalities defining the original generalized polyhedron.
% This result is a natural extension of the corresponding result for closed polyhedra~\cite{Schrijver1998}.
\begin{lemma}
\label{lemma:Projection}
Let $n \in \PositiveIntegers{}$ and let $P \in \mathbb{R}^{n}$ be a rational generalized polyhedron defined by inequalities with coefficients of size at most $\MaximumEntrySize$.
Then, the projection of $P$ onto any subspace is a rational generalized polyhedron defined by inequalities with coefficients of size polynomial in $n$ and $\MaximumEntrySize$.
\end{lemma}

\begin{MyProof}
The claim holds trivially when $n = 1$.
Thus, we assume $n \ge 2$.
Let $n_x, n_y \in \PositiveIntegers{}$ and $m_{\ge}, m_{>} \in \NonNegativeIntegers{}$, where $n = n_x + n_y$.
Let $A \in \mathbb{Q}^{n_x \times m_{\ge}}$, $B \in \mathbb{Q}^{n_y \times m_{\ge}}$, $C \in \mathbb{Q}^{n_x \times m_{>}}$, $D \in \mathbb{Q}^{n_y \times m_{>}}$, $a \in \mathbb{Q}^{m_{\ge}}$, and $c \in \mathbb{Q}^{m_{>}}$ be such that $P = \{ (x, y) \in \mathbb{R}^{n_x + n_y} : A x + B y \ge a, C x + D y > c \}$, where the size of each coefficient is at most $\MaximumEntrySize$.
We show that the claim holds for $\mathrm{proj}_x(P)$.
Define the projection cone 
$$
C_P = \{ (u, v) \in \mathbb{R}^{m_{\ge}} \times \mathbb{R}^{m_{>}}: B^{\top} u + D^{\top} v = \ZeroVector, u \ge \ZeroVector, v \ge \ZeroVector \},
$$
which is also used in the context of the projection of the closed polyhedron~\cite{conforti2014integer}.
Let $\{ (d_i, e_i) \in \mathbb{R}^{m_{\ge}} \times \mathbb{R}^{m_{>}} : i \in I\}$ be the extreme rays of $C_P$ and define $I_1 = \{ i \in I : e_i = \ZeroVector \}$ and $I_2 = \{ i \in I : e_i \ne \ZeroVector \}$.
The claim holds trivially if $P = \emptyset$.
Moreover, $C_P = \{ \ZeroVector \}$ implies $\mathrm{proj}_x(P) = \mathbb{R}^{n_x}$.
Thus, we assume $P \ne \emptyset$ and $C_P \ne \{ \ZeroVector \}$.

For each $\epsilon > 0$, define the perturbed set: $P_\epsilon = \{ (x, y) \in \mathbb{R}^{n_x + n_y} : A x + B y \ge a, C x + D y \ge c + \epsilon \OneVector \}$.
% Since $P_\epsilon$ is a closed polyhedron for any $\epsilon > 0$, 
A standard result in polyhedral theory (see~\cite{conforti2014integer}) gives:
$$
\mathrm{proj}_x(P_\epsilon) = \left\{
x \in \mathbb{R}^{n_x}
:
\begin{array}{ll}
d_i^{\top} A x \ge d_i^{\top} a,
&
\forall i \in I_1, \\
d_i^{\top} A x + e_i^{\top} C x \ge d_i^{\top} a + e_i^{\top} (c + \epsilon \OneVector),
&
\forall i \in I_2, \\
\end{array}
\right\}.
$$
We have
\begin{align}
\mathrm{proj}_x(P) 
&=
\mathrm{proj}_x\left(\bigcup_{\epsilon > 0} P_\epsilon\right)
\notag
\\
&= 
\bigcup_{\epsilon > 0} \mathrm{proj}_x(P_\epsilon)
\notag
\\
&=
\bigcup_{\epsilon > 0}
\left\{
x \in \mathbb{R}^{n_x}
:
\begin{array}{ll}
d_i^{\top} A x \ge d_i^{\top} a,
&
\forall i \in I_1, \\
d_i^{\top} A x + e_i^{\top} C x \ge d_i^{\top} a + e_i^{\top} (c + \epsilon \OneVector),
&
\forall i \in I_2, \\
\end{array}
\right\}
\notag
\\
&=
\left\{
x \in \mathbb{R}^{n_x}
:
\begin{array}{ll}
d_i^{\top} A x \ge d_i^{\top} a,
&
\forall i \in I_1, \\
d_i^{\top} A x + e_i^{\top} C x > d_i^{\top} a + e_i^{\top} c,
&
\forall i \in I_2, \\
\end{array}
\right\},
\label{eq:LemmaProjectionFormula}
\end{align}
which is a rational generalized polyhedron.
% Exercise 4.40 in Bertsimas.
To bound the size of the coefficients in the inequalities in~\eqref{eq:LemmaProjectionFormula}, we will bound the number and the size of the nonzero entries of each extreme ray of $C_P$.
The set of the extreme rays of $C_P$ corresponds precisely to the set of extreme points of the set $D_P$~\cite{bertsimas1997introduction}, where
$$
D_P = \{ (u, v) \in \mathbb{R}^{m_{\ge}} \times \mathbb{R}^{m_{>}}: B^{\top} u + D^{\top} v = \ZeroVector, u \ge \ZeroVector, v \ge \ZeroVector, \OneVector^{\top} u + \OneVector^{\top} v = 1 \}.
$$
Let $(u, v) \in D_P$ be an extreme point of $D_P$.
Then, there exist $m_{\ge} + m_{>}$ linearly independent constraints active at $(u, v)$.
Since $D_P$ contains $n_y + 1$ equality constraints and $m_{\ge} + m_{>}$ inequality constraints (namely, $u \ge \ZeroVector$ and $v \ge \ZeroVector$), at least $m_{\ge} + m_{>} - n_y - 1$ of the inequalities must be active at $(u, v)$.
That is, at most $n_y + 1$ inequalities can be inactive.
Therefore, any extreme point of $D_P$ has at most $n_y + 1$ nonzero entries.
Now, let $(u', v')$ be the subvector of nonzero entries, and let $B', D', b'$ denote the corresponding coefficients and RHS of linearly independent active constraints, so that:
$$
(B' \ D') \begin{pmatrix} u' \\ v' \end{pmatrix} = b'.
$$
Note that $(u', v')$ is the unique solution to this system of linear equations.
Thus, it is rational of size polynomially bounded in the number of variables $n_y + 1$ and the maximum size of the entries, which is at most $\MaximumEntrySize$~\cite{Schrijver1998}.
Thus, the entries in the inequalities in~\eqref{eq:LemmaProjectionFormula} defining $\mathrm{proj}_x(P)$ have size polynomial in $n$ and $\MaximumEntrySize$.
\MyQED
\end{MyProof}

We state a simple result concerning the complement of a rational generalized polyhedron.

\begin{lemma}
\label{lemma:ComplementAsUnion}
Let $P$ be a rational generalized polyhedron defined by inequalities with coefficients of size at most $\MaximumEntrySize$.
Then, the complement of $P$ is a disjoint union of rational generalized polyhedra defined by inequalities with coefficients of size at most $\MaximumEntrySize$.
\end{lemma}

\begin{MyProof}
Let $P = \{ x \in \mathbb{R}^n : A x \ge b, C x > d \}$, where $A \in \mathbb{Q}^{m \times n}$, $b \in \mathbb{Q}^{m}$, $C \in \mathbb{Q}^{m' \times n}$, $d \in \mathbb{Q}^{m'}$.
We denote the $j$-th row of $A$ and $C$ by $a_j$ and $c_j$, respectively.
Then, 
% the complement of $P$ is
\begin{align*}
\mathbb{R}^{n} \setminus P
&=
\{ x \in \mathbb{R}^{n} : a_1^{\top} x < b_1 \}
\cup
\{ x \in \mathbb{R}^{n} : a_1^{\top} x \ge b_1, a_2^{\top} x < b_2 \}
\\
&\qquad
\cup \cdots \cup \{ x \in \mathbb{R}^{n} : a_1^{\top} x \ge b_1, \ldots, a_{m - 1}^{\top} x \ge b_{m - 1}, a_{m}^{\top} x < b_{m} \}
\\
&\qquad
\cup \{ x \in \mathbb{R}^{n} : A x \ge b, c_1^{\top} x \le d_1 \}
\cup \{ x \in \mathbb{R}^{n} : A x \ge b, c_1^{\top} x > d_1, c_2^{\top} \le d_2 \}
\\
&\qquad
\cup \cdots \cup \{ x \in \mathbb{R}^{n} : A x \ge b, c_1^{\top} x > d_1, \ldots, c_{m' - 1}^{\top} x > d_{m' - 1}, c_{m'}^{\top} x \le d_{m'} \},
\end{align*}
from which the claim follows.
% which is a disjoint union of rational generalized polyhedra defined by inequalities with coefficients of size at most $\MaximumEntrySize$.
\MyQED
\end{MyProof}

We say a function $\theta: \mathbb{R}^{n} \rightarrow \mathbb{R} \cup \{\pm\infty\}$ is a \emph{rational piecewise linear function} if and only if $\theta$ can be written as
\begin{equation*}
\theta(x) = \begin{cases}
c_l^{\top} x, & \text{ if $x \in P_l$ for some $l \in \mathcal{L}$}, \\
\infty, & \text{ if $x \in P_l$ for some $l \in \mathcal{L}^+$}, \\
-\infty, & \text{ if $x \in P_l$ for some $l \in \mathcal{L}^-$},
\end{cases}
\end{equation*}
where $\mathcal{L}$, $\mathcal{L}^+$ and $\mathcal{L}^-$ are finite, disjoint index sets, $c_l \in \mathbb{Q}^n$ for each $l \in \mathcal{L}$, and $\{ P_l : l \in \mathcal{L} \cup \mathcal{L}^+ \cup \mathcal{L}^- \}$ is a partition of $\mathbb{R}^{n}$ into rational generalized polyhedra.
We say that $\theta$ is defined by inequalities of size at most $\MaximumEntrySize$ if the size of each entry in $c_l$, $l \in \mathcal{L}$, is at most $\MaximumEntrySize$, and each of $P_l$ is defined by inequalities with coefficients of size at most $\MaximumEntrySize$.
% In this case, the epigraph of $\theta$ is the union of rational generalized polyhedra defined by inequalities with coefficients of size at most $\MaximumEntrySize$.
% \NagisaSideComment{\textcolor{blue!40}{The epigraph may be defined by inequalities of smaller sizes than the original piecewise linear function.}}

\begin{lemma}
\label{lemma:MinimumOfPiecewiseLinearFunctions}
Let $m \in \PositiveIntegers{}$, and let $\{ \theta_i : i = 1, \ldots, m \}$ be a collection of rational piecewise linear functions defined by inequalities with coefficients of size at most $\MaximumEntrySize$.
Then, $\theta(z) = \min_{i = 1, \ldots, m } \theta_i(z)$ is a rational piecewise linear function defined by inequalities with coefficients of size at most $\MaximumEntrySize + 1$.
\end{lemma}

\begin{MyProof}
For each $i = 1, \ldots, m$, write
\begin{equation}
\theta_i(x) = \begin{cases}
c_{i l_i}^{\top} x, & \text{ if $x \in P_{l_i}$ for some $l_i \in \mathcal{L}_i$}, \\
\infty, & \text{ if $x \in P_{l_i}$ for some $l_i \in \mathcal{L}^+_i$}, \\
-\infty, & \text{ if $x \in P_{l_i}$ for some $l_i\in \mathcal{L}^-_i$},
\end{cases}
\notag
\end{equation}
where $\mathcal{L}_i$, $\mathcal{L}^+_i$ and $\mathcal{L}^-_i$ are finite, disjoint index sets, $c_{i l_i} \in \mathbb{Q}^n$ for each $l_i \in \mathcal{L}_i$ and $\{ P_{i l_i} : l_i \in \mathcal{L}_i \cup \mathcal{L}^+_i \cup \mathcal{L}^-_i \}$ is a partition of $\mathbb{R}^n$ into rational generalized polyhedra.
In particular, we can assume that, for all $i = 1, \ldots, m$, and $l \in \mathcal{L}_i$, the sizes of entries in $c_{i l_i}$ are at most $\MaximumEntrySize$, and each $P_{i l_i}$ is defined by rational inequalities with coefficients of size at most $\MaximumEntrySize$.
Let $\bar{\mathcal{L}}_i = \mathcal{L}_i \cup \mathcal{L}^+_i \cup \mathcal{L}^-_i$ for each $i = 1, \ldots, m$, and let $\bar{\mathcal{L}} = \bar{\mathcal{L}}_1 \times \cdots \times \bar{\mathcal{L}}_m$.
For each $l = (l_1, \ldots, l_m) \in \bar{\mathcal{L}}$, define
$$
R_l = P_{1 l_1} \cap \cdots \cap P_{m l_m}.
$$
Then, $\{ R_l : l \in \bar{\mathcal{L}} \}$ is a disjoint collection of rational generalized polyhedra that covers $\mathbb{R}^n$, with each $R_l$ defined by inequalities with coefficients of size at most $\MaximumEntrySize$.
Now we further decompose $R_l$ into disjoint regions $\{ R_l^i : i = 1, \ldots, m \}$ as follows: if $R_l = \emptyset$, define $R_l^i = \emptyset$ for all $i$.
Otherwise, for each $i = 1, \ldots, m$, define
\begin{align*}
R_l^i 
&= \{ x \in R_l : \theta_i(x) = \min_{j = 1, \ldots, m} \theta_j(x), \ \theta_j(x) > \theta_i(x), \ \forall j < i \} \\
&= \{ x \in R_l : \theta_j(x) > \theta_i(x), \ \forall j < i, \  \theta_j(x) \ge \theta_i(x), \ \forall j > i \}.
\end{align*}
That is, $R_l^i$ consists of the points in $R_l$ where $\theta_i$ achieves the minimal value among $\theta_1, \ldots, \theta_m$, and ties are broken in favor of the smallest index.
Each of $R_l^i$ is defined by inequalities with coefficients of size at most $\MaximumEntrySize + 1$.
For example, in the case where all $\theta_i$ take finite values over $R_l$, i.e., $l_i \in \mathcal{L}_i$ for all $i$, we have:
$$
R_l^i = \{ x \in R_l : c_{j l_{j}}^{\top} x > c_{i l_i}^{\top} x, \ \forall j < i, \ \ c_{j l_{j}}^{\top} x \ge c_{i l_i}^{\top} x, \ \forall j > i \},
$$
whose coefficients are of size at most $\MaximumEntrySize + 1$.
The cases where $\theta_i = \pm\infty$ are handled similarly.

Now, for any $x \in \mathbb{R}^n$, there exists a unique $l \in \mathcal{L}$ and unique $i \in \{1, \ldots, m\}$ such that $x \in R_l^i$, and in this region we have $\theta(x) = \theta_i(x)$, which is a linear function on $R_l^i$.
Thus, $\theta$ is a rational piecewise linear function, linear over each generalized polyhedron in $\{R_l^i : l \in \bar{\mathcal{L}}, i = 1, \ldots, m \text{ such that } R_l^i \ne \emptyset\}$, each of which is described by inequalities with coefficients of size at most $\MaximumEntrySize + 1$.
\MyQED
\end{MyProof}

A rational piecewise linear function can arise as the value function of an LP, or more generally, optimization of a linear function over a union of rational generalized polyhedra.
\begin{lemma}
\label{lemma:ValueFunctionOfLinearOptimizationOverTheUnionOfGeneralizedPolyhedra}
Let $n_x, n_y \in \PositiveIntegers{}$ and define $n = n_x + n_y$.
% ObjTruncation
Let $c \in \mathbb{Q}^{n_x}$ be a rational vector with entries of size at most $\MaximumEntrySize$, and let $\{ P_l \subset \mathbb{R}^n : l \in \mathcal{L} \}$ be a finite collection of rational generalized polyhedra, each defined by inequalities with coefficients of size at most $\MaximumEntrySize$.
% ObjTruncation
Then, the function $v(y) = \inf_{x} \{ c^{\top} x : (x, y) \in \cup_{l \in \mathcal{L}} P_l \}$ is a rational piecewise linear function defined by inequalities with coefficients of size polynomial in $\NVariables$ and $\MaximumEntrySize$.
% at most $\ValueFunctionOfUnionLPPolynomial(\MaximumEntrySize, n)$, for some polynomial $\ValueFunctionOfUnionLPPolynomial$.
\end{lemma}

\begin{MyProof}
We first consider the case where $|\mathcal{L}| = 1$, i.e., the feasible set is a single rational generalized polyhedron $P \subset \mathbb{R}^n$.
If $P = \emptyset$, then $v(y) = \infty$, and the claim holds trivially.
Assume instead that $P \ne \emptyset$.
Let $\bar{P} = \mathrm{cl}(P)$ be the closure of $P$, which is a rational polyhedron (Lemma~\ref{lemma:estimate-of-weak-and-strong-linear-inequality}~\ref{lemma:estimate-of-weak-and-strong-linear-inequality:one}).
Define the function 
% ObjTruncation
$$
\bar{v}(y) = \inf_{x} \{c^{\top} x : (x, y) \in \bar{P} \}.
$$
Being the value function of an LP, $\bar{v}(y)$ is a rational piecewise linear function defined by inequalities with coefficients of size polynomial in $n$ and $\MaximumEntrySize$~(see, e.g., \cite{ward1990approaches,Schrijver1998}).
Now, observe that 
$$
v(y) = \begin{cases}
\bar{v}(y), & \text{ if $y \in \mathrm{proj}_y(P),$} \\
\infty, & \text{ otherwise.}
\end{cases}
$$
The claim follows by Lemmata~\ref{lemma:Projection} and~\ref{lemma:ComplementAsUnion}.

Now consider the general case where $|\mathcal{L}| > 1$.
For each $l \in \mathcal{L}$, define the function 
% ObjTruncation
$$
v_l(y) = \inf_x \{ c^{\top} x : (x, y) \in P_l \}.
$$
As shown above, for each $l \in \mathcal{L}$, $v_l$ is a rational piecewise linear function defined by inequalities with coefficients of size polynomial in $\MaximumEntrySize$ and $n$.
Since $v(y) = \min_{l \in \mathcal{L}} v_l(y)$, the result follows from Lemma~\ref{lemma:MinimumOfPiecewiseLinearFunctions}.
\MyQED
\end{MyProof}

\subsection{Feasible Set of $k$-level LP}
\label{sec:FeasibleSetOfMultilevelLP}

In this section, we establish an estimate on the size of a feasible and optimal solution to a $k$-level LP instance for any fixed $k \ge 1$.
To this end, we show that the feasible set of a $k$-level LP instance is a union of rational generalized polyhedra.

We first extend the value function approach used for bilevel programming to $k$-level LP.
\begin{lemma}
\label{lemma:ValueFunctionReformulation}
For any $l = 1, \ldots, k$ and $x_1, \ldots, x_{l - 1}$, it holds that
% ObjTruncation
\begin{align}
&
\mathcal{F}\EqrefGeneralKLP{k}{l}{x_1, \ldots, x_{l - 1}, A, b, c} =
\notag
\\
&
\qquad
\left\{ 
\begin{pmatrix}
x_l \\ \vdots \\ x_k
\end{pmatrix}
:
\begin{array}{ll}
\displaystyle
\sum_{i = 1}^k A_{l' i} x_i \ge b_{l'},
&
\forall l' = l, \ldots, k,
\\
\displaystyle
\sum_{i = l'}^k c_{l' i}^{\top} x_i \le v\EqrefGeneralKLP{k}{l'}{x_1, \ldots, x_{l' - 1}, A, b, c}, 
&
\forall l' = l + 1, \ldots, k
\end{array}
\right\}.
\label{eq:lemma:ValueFunctionReformulation}
% \label{eq:LLevelValueFunctionReformulation}
% \tag{$\mathrm{P}_{l}(x_{l + 1}, \ldots, x_k)$}
\end{align}
\end{lemma}

\begin{MyProof}
We prove by induction, starting from $l = k$.
The base case $l = k$ holds by definition.
Now, suppose the claim holds for some $l + 1$.
Pick any values of $x_1, \ldots, x_{l - 1}$.
A solution $(x_l, \ldots, x_k)$ is feasible to the $l$-th player's problem if and only if it satisfies the linear constraints in the $l$-th player's problem~\eqref{eq:lemma:ValueFunctionReformulation:one} and $(x_{l + 1}, \ldots, x_k)$ is optimal to the $(l+1)$-th player's problem given $x_1, \ldots, x_{l - 1}$ and $x_l$.
The latter condition holds if and only if $(x_{l + 1}, \ldots, x_k)$ is feasible to the $(l+1)$-th player's problem (by induction hypothesis, it is equivalent to satisfying~\eqref{eq:lemma:ValueFunctionReformulation:two} and \eqref{eq:lemma:ValueFunctionReformulation:three}) and it attains the optimal objective value~\eqref{eq:lemma:ValueFunctionReformulation:four}.
Thus, $(x_l, \ldots, x_k)$ is feasible to the $l$-th player's problem if and only if
% ObjTruncation
\begin{align}
\sum_{i = 1}^k A_{l i} x_i &\ge b_l,
\label{eq:lemma:ValueFunctionReformulation:one}
\\
\sum_{i = 1}^k A_{l' i} x_i &\ge b_{l'},
&
\forall l' = l + 1, \ldots, k,
\label{eq:lemma:ValueFunctionReformulation:two}
\\
\sum_{i = l'}^k c_{l' i}^{\top} x_i &\le v\EqrefGeneralKLP{k}{l'}{x_1, \ldots, x_{l' - 1}, A, b, c}, 
&
\forall l' = l + 2, \ldots, k
\label{eq:lemma:ValueFunctionReformulation:three}
\\
\sum_{i = l + 1}^k c_{l + 1 \, i \,}^{\top} x_i &\le v\EqrefGeneralKLP{k}{l + 1}{x_1, \ldots, x_l, A, b, c}, 
\label{eq:lemma:ValueFunctionReformulation:four}
% \label{eq:LLevelValueFunctionReformulation}
% \tag{$\mathrm{P}_{l}(x_{l + 1}, \ldots, x_k)$}
\end{align}
implying the claim holds for $l$.
\MyQED
\end{MyProof}

\begin{lemma}
\label{lemma:PiecewiseLinearityOfValueFunction}
Let $\MaximumEntrySize$ be the maximum size of entries in $A$, $b$ and $c$.
For each $l = 2, \ldots, k$, the value function of the $l$-th player's problem, $v\EqrefGeneralKLP{k}{l}{x_1, \ldots, x_{l - 1}, A, b, c}$, is a rational piecewise linear function in $x_1, \ldots, x_{l - 1}$ defined by inequalities with coefficients of size polynomial in $\NVariables$ and $\MaximumEntrySize$.
\end{lemma}

\begin{MyProof}
We prove by induction, starting from $l = k$.
The base case $l = k$ follows from Lemma~\ref{lemma:ValueFunctionOfLinearOptimizationOverTheUnionOfGeneralizedPolyhedra}.
Suppose the claim holds for some $l + 1$.
Then, by Lemma~\ref{lemma:ValueFunctionReformulation}, $\mathcal{F}\EqrefGeneralKLP{k}{l}{x_1, \ldots, x_{l - 1}, A, b, c}$ is given by \eqref{eq:lemma:ValueFunctionReformulation}.
By induction hypothesis, for each $l' = l + 1, \ldots, k$, $v\EqrefGeneralKLP{k}{l'}{x_1, \ldots, x_{l'-1}, A, b, c}$ is a rational piecewise linear function in $x_1, \ldots, x_{l'-1}$ defined by inequalities with coefficients of polynomial size.
Thus, $\mathcal{F}\EqrefGeneralKLP{k}{l}{x_1, \ldots, x_{l - 1}, A, b, c}$ (more precisely, the graph of the set-valued mapping $(x_1, \ldots, x_{l - 1}) \mapsto \mathcal{F}\EqrefGeneralKLP{k}{l}{x_1, \ldots, x_{l - 1}, A, b, c}$) is a union of generalized polyhedra.
Therefore, by Lemma~\ref{lemma:ValueFunctionOfLinearOptimizationOverTheUnionOfGeneralizedPolyhedra}, $v\EqrefGeneralKLP{k}{l}{x_1, \ldots, x_{l - 1}, A, b, c}$ is a rational piecewise linear function in $x_1, \ldots, x_{l-1}$ defined by inequalities with coefficients of size polynomial in $\NVariables$ and $\MaximumEntrySize$.
\MyQED
\end{MyProof}

The following lemma gives an estimate of the size of a solution to a $k$-level LP instance, which is crucial to studying the complexity of $k$-level LP.

\begin{lemma}
\label{lemma:SolutionSizeEstimate}
% \mbox{}
Let $(A, b, c)$ be a rational $k$-level LP instance.
\begin{enumerate}
\item
\label{lemma:SolutionSizeEstimate:One}
Its feasible set $\mathcal{F}\eqref{LabelGeneralKLP}$ is the union of generalized polyhedra defined by inequalities with coefficients of polynomial size.
\item
\label{lemma:SolutionSizeEstimate:Two}
If it is feasible, it has a rational feasible solution of size polynomial in the input size.
\item
If it has an optimal solution, it has a rational optimal solution of size polynomial in the input size.
\item
If its optimal objective value is finite, it is a rational number of size polynomial in the input size.
\end{enumerate}
\end{lemma}

\begin{MyProof}
Invoke Lemmata~\ref{lemma:estimate-of-weak-and-strong-linear-inequality},~\ref{lemma:ValueFunctionReformulation} and~\ref{lemma:PiecewiseLinearityOfValueFunction}.
\MyQED
\end{MyProof}

Lemma~\ref{lemma:SolutionSizeEstimate} states that if the instance has an optimal solution, there exists a rational optimal solution of polynomial size.
This stands in sharp contrast to the classical three-player game of \citet{Nash1951}, whose unique mixed Nash equilibrium is irrational.
Moreover, Lemma~\ref{lemma:SolutionSizeEstimate} implies that, if the instance is unbounded, then there exist a rational point $(x_1, \ldots, x_k)$ and a rational direction vector $(d_1, \ldots, d_k)$, both of size polynomial in the input size, such that $(x_1, \ldots, x_k) + t (d_1, \ldots, d_k)$ is feasible for all $t \ge 0$ and the objective value decreases as $t$ increases.
Clearly, in general, no polynomial-time verifier can certify the correctness of such a witness $(x_1, \ldots, x_k)$ and $(d_1, \ldots, d_k)$, unless the polynomial hierarchy collapses to its first level.

\subsection{Main Analysis}
\label{sec:ProofOfInclusionInSigmaPKMinusOne:MainAnalysis}

In this section, we establish the inclusion of \DecisionProblemVal{k} and \DecisionProblemUnb{k} in \ClassSigmaP{k - 1} for any fixed $k \ge 1$.
To this end, we consider the following two decision problems.
\begin{align}
&
\parbox{0.75\textwidth}{Given a rational $k$-level LP instance $(A, b, c)$ and a solution $(x_1, \ldots, x_k)$, is it optimal, i.e., $(x_1, \ldots, x_k) \in \OptSolSet\eqref{LabelGeneralKLP}$?}
\tag{\EqTagDecisionProblemOnKLP{k}{RECOG}}
\label{LabelDecisionProblemRecognize}
\\
&
\parbox{0.75\textwidth}{Given a rational $k$-level LP instance $(A, b, c)$, is it feasible?}
\tag{\EqTagDecisionProblemOnKLP{k}{FEAS}}
\label{LabelDecisionProblemFeas}
\end{align}

We investigate each problem in order.

\begin{lemma}
\label{lemma:ComplexityOfRecognition}
The decision problem \DecisionProblemRecognize{k} is in \ClassPiP{k - 1}.
\end{lemma}
\begin{MyProof}
We prove by induction.
The base case $k = 1$ follows from polynomial solvability of LP~\cite{Schrijver1998}.
Suppose the claim holds for $k - 1$.
We show that \DecisionProblemModifierComp{\DecisionProblemRecognize{k}}, the complement of \DecisionProblemRecognize{k}, is in \ClassSigmaP{k - 1}.
Let $(A, b, c)$ be a rational $k$-level LP instance, and let $(\bar{x}_1, \ldots, \bar{x}_k)$ be a candidate solution.
This point is not an optimal solution if and only if one of the following conditions holds.
\begin{enumerate}
\item
Solution $(\bar{x}_1, \ldots, \bar{x}_k)$ violates the linear constraints in the first player's problem;
\item 
Solution $(\bar{x}_2, \ldots, \bar{x}_k)$ is not optimal to the second player given $\bar{x}_1$;
\item 
There exists a solution $(\hat{x}_1, \ldots, \hat{x}_k)$ feasible to the first player's problem with a strictly smaller objective value than $(\bar{x}_1, \ldots, \bar{x}_k)$.
\end{enumerate}
Note that in case (iii), there exists a rational $(\hat{x}_1, \ldots, \hat{x}_{k})$ of polynomial size.
This follows from the fact that the set
$$
P = \left\{ 
\begin{pmatrix}
x_1 \\ \ldots \\ x_k
\end{pmatrix}
\in \mathcal{F}\eqref{LabelGeneralKLP} : \sum_{i = 1}^k c_{1 i}^{\top} x_i < \sum_{i = 1}^k c_{1 i}^{\top} \bar{x}_i \right\}
$$
is a finite union of rational generalized polyhedra, each defined by inequalities with coefficients of size polynomial in the input size (including the size of $(\bar{x}_1, \ldots, \bar{x}_k)$).
Thus, if $P$ is nonempty, then it contains a rational point of polynomial size.
Based on this, we construct a nondeterministic Turing machine (NDTM) $N$ with oracle access to \DecisionProblemRecognize{(k - 1)} that decides \DecisionProblemModifierComp{\DecisionProblemRecognize{k}} as follows.
\begin{itemize}
\item
Nondeterministically guess which of the three cases applies.
\item
If case (i) is guessed, verify whether the candidate violates any linear constraints in the first player's problem. Accept if it does, and reject otherwise.
\item
If case (ii) is guessed, query the oracle to check whether $(\bar{x}_2, \ldots, \bar{x}_k)$ is optimal for the second player's problem given $\bar{x}_1$. Accept if the oracle returns ``no''; otherwise, reject.
\item
If case (iii) is guessed, nondeterministically guess a rational solution $(\hat{x}_1, \ldots, \hat{x}_k)$ of polynomial size, verify that it has a strictly smaller objective value and satisfies the linear constraints in the first player's problem, and use the oracle to check whether $(\hat{x}_2, \ldots, \hat{x}_k)$ is optimal for the second player's problem given $\hat{x}_1$. Accept if all checks pass, and reject otherwise.
\end{itemize}
The instance is a YES instance of \DecisionProblemModifierComp{\DecisionProblemRecognize{k}} if and only if there exist nondeterministic choices that make $N$ accept, implying \DecisionProblemRecognize{k} is in \ClassPiP{k-1}.
\MyQED
\end{MyProof}

\begin{lemma}
\label{lemma:ComplexityOfFeas}
The decision problem \DecisionProblemFeas{k} is in \ClassSigmaP{k - 1}.
\end{lemma}

\begin{MyProof}
The claim holds for $k = 1$ from the polynomial solvability of LP~\cite{Schrijver1998}.
In the following, assume $k \ge 2$.
Let $(A, b, c)$ be a rational $k$-level LP instance.
A solution $(x_1, \ldots, x_k)$ is feasible if and only if the following two conditions are satisfied:
\begin{enumerate}
\item
Solution $(x_1, \ldots, x_k)$ satisfies the linear constraints in the first player's problem;
\item 
Solution $(x_2, \ldots, x_k)$ is optimal for the second player's problem given $x_1$.
\end{enumerate}
By Lemma~\ref{lemma:SolutionSizeEstimate}, if the input instance is feasible, it has a rational feasible solution $(x_1, \ldots, x_k)$ of polynomial size.
Based on this, we construct an NDTM $M$ with oracle access to \DecisionProblemRecognize{(k - 1)} that decides \DecisionProblemFeas{k} as follows.
\begin{itemize}
\item
Nondeterministically guess a rational feasible solution $(x_1, \ldots, x_k)$ of polynomial size.
\item
Verify that it satisfies the linear constraints in the first player's problem, and use the oracle to check whether $(x_2, \ldots, x_k)$ is optimal for the second player's problem given $x_1$. Accept if both checks pass.
\end{itemize}
The instance is a YES instance of \DecisionProblemFeas{k} if and only if there exist nondeterministic choices under which $M$ accepts.
Now the claim follows in light of Lemma~\ref{lemma:ComplexityOfRecognition}.
\MyQED
\end{MyProof}

\begin{lemma}
\label{lemma:ValAndDecisionProblemUnbsinSigmaPKminusOne}
\mbox{}
\begin{enumerate}
\item 
\label{lemma:ValAndDecisionProblemUnbsinSigmaPKminusOne:One}
The decision problem \DecisionProblemVal{k} is in \ClassSigmaP{k - 1}.
\item
\label{lemma:ValAndDecisionProblemUnbsinSigmaPKminusOne:Two}
The decision problem \DecisionProblemUnb{k} is in \ClassSigmaP{k - 1}.
\end{enumerate}
\end{lemma}

\begin{MyProof}
\mbox{} % // TODO
\begin{enumerate}
\item 
Let $(A, b, c)$ be a rational instance, and let $t$ be a rational number.
We have $v\eqref{LabelGeneralKLP} \le t$ if and only if the following instance is feasible:
\begin{equation*}
\inf_{x_1, \ldots, x_k}
\left\{ 
\sum_{i = 1}^k c_{1 i}^{\top} x_i 
: 
% \begin{array}{l}
% \displaystyle
\sum_{i = 1}^k A_{1 i} x_i \ge b_1, 
\sum_{i = 1}^k c_{1 i}^{\top} x_i \le t,
% \\
% (x_2, \ldots, x_k)
\begin{pmatrix}
x_2 \\ \vdots \\ x_k
\end{pmatrix}
\in 
\mathcal{S}\EqrefGeneralKLP{k}{2}{x_1, A, b, c}
% \end{array}
\right\},
\end{equation*}
Thus, \DecisionProblemVal{k} $\le_l$ \DecisionProblemFeas{k}.
The claim follows from Lemma~\ref{lemma:ComplexityOfFeas}.
\item
Let $(A, b, c)$ be a rational instance.
Let $\phi$ be the polynomial that bounds the size of the optimal objective value if it is finite (Lemma~\ref{lemma:SolutionSizeEstimate}).
Instance $(A, b, c)$ is unbounded if and only if $v\eqref{LabelGeneralKLP} \le -2^{\phi(l)} - 1$, where $l$ is the input size.
Thus, \DecisionProblemUnb{k} $\le_l$ \DecisionProblemVal{k}.
The claim follows from assertion~\ref{lemma:ValAndDecisionProblemUnbsinSigmaPKminusOne:One}.
\MyQED
\end{enumerate}
\end{MyProof}

\section{Extensions and Limitations}
\label{sec:extensions_and_limitations}

The analysis in Section~\ref{sec:ProofOfMembership} readily extends to the case in which $\EqrefGeneralKLP{k}{1}{A, b, c}$ includes strict inequalities of the form $\sum_{i=1}^k A'_{l i} x_i > b'_l$ for each $l$. We refer to this variant as a generalized $k$-level LP. In particular, the decision version of the generalized $k$-level LP belongs to $\ClassSigmaP{k-1}$. Moreover, the version of this decision problem that asks whether there exists a feasible solution whose objective value is strictly better than a given threshold (i.e., using strict rather than weak inequalities) also lies in $\ClassSigmaP{k-1}$. It is well known that the inclusion of strict inequalities does not affect the polynomial-time solvability of (single-level) LPs~\cite{Schrijver1998}. The argument above shows that, similarly, the presence of strict inequalities does not affect the complexity of multilevel LPs established in this paper.
% For clarity and simplicity, we use the standard notation $\EqrefGeneralKLP{k}{1}{A, b, c}$ throughout the main text.

A more practically motivated extension is to consider the set of $\epsilon$-optimal solutions.
In multilevel programming, the exact optimality of the followers is often relaxed, e.g., to ensure solvability of the leader’s problem~\cite{loridan1996weak}.
For an optimization instance $E$, define
$$
\OptSolSet_{\epsilon}(E) = \{ x : x \in \FeasSolSet(E) \text{ s.t. its objective value is less than or equal to $\OptValFunc(E) + \epsilon$} \}.
$$
Suppose that in \eqref{LabelGeneralKLP}--\eqref{eq:1Level}, we replace $\OptSolSet$ by $\OptSolSet_{\epsilon}$.
Then a variant of Lemma~\ref{lemma:ValueFunctionReformulation} holds: the feasible set of the $l$-th player $\mathcal{F}\EqrefGeneralKLP{k}{l}{x_1, \ldots, x_{l - 1}, A, b, c}$ is given by
% ObjTruncation
\begin{align*}
% &
%  =
% \notag
% \\
% &
% \qquad
\left\{ 
\begin{pmatrix}
x_l \\ \vdots \\ x_k
\end{pmatrix}
:
\begin{array}{ll}
\displaystyle
\sum_{i = 1}^k A_{l' i} x_i \ge b_{l'},
&
\forall l' = l, \ldots, k,
\\
\displaystyle
\sum_{i = l'}^k c_{l' i}^{\top} x_i \le v\EqrefGeneralKLP{k}{l'}{x_1, \ldots, x_{l' - 1}, A, b, c} + \epsilon, 
&
\forall l' = l + 1, \ldots, k
\end{array}
\right\}.
\end{align*}
Consequently, all arguments in our analysis carry over mutatis mutandis, and the decision version of a $k$-level LP with $\epsilon$-optimal solution sets lies in $\ClassSigmaP{k-1}$.
A similar extension is possible if we define the set of $\epsilon$-optimal solutions in the relative sense rather than in the absolute sense.

From the arguments of \citet{Jeroslow1985}, it follows that multilevel pure binary LP with an arbitrary number of levels is \ClassPSPACE{}-complete.
For multilevel (continuous) LP, however, it remains open whether the decision problem with unbounded $k$ lies in \ClassPSPACE{} (\ClassPSPACE{}-hardness is immediate).
The difficulty stems from the fact that our analysis, including Lemma~\ref{lemma:SolutionSizeEstimate}, does not extend to the setting in which $k$ is unbounded.
From a theoretical standpoint, determining the complexity of the decision version of multilevel LP with an arbitrary number of levels is an important open problem.

On the positive side, our analysis naturally extends to show that the decision version of $k$-level mixed binary LP is $\ClassSigmaP{k}$-complete for every $k \ge 1$, and remains so under conditions~\ref{Condition1}--\ref{Condition3}. 
Membership in $\ClassSigmaP{k}$ follows from the observation that any such instance can be reformulated as a $(k+1)$-level continuous LP. 
This reformulation is achieved by applying the technique of \citet{audet1997links}, which transforms a binary LP into a bilevel LP. 
Hardness follows from the fact that the decision version of the $k$-level pure binary LP is $\ClassSigmaP{k}$-hard, even under Conditions~\ref{Condition1}--\ref{Condition3}, as shown by \citet{Jeroslow1985}.
% Likewise, any $k$-level mixed integer LP instance in which all variables have finite (not necessarily polynomial) lower and upper bounds can be reformulated as a $k$-level mixed binary LP instance using binary expansion.
% Consequently, for each $k$, $k$-level mixed integer LP with finite variable bounds is $\ClassSigmaP{k}$-complete.

\section{Conclusions}
\label{sec:Conclusions}

This work has provided a comprehensive analysis of the computational complexity of decision problems associated with $k$-level LP.
In particular, we proved that the decision version of $k$-level LP, determining if the optimal objective value is equal to or better than a given threshold, belongs to $\ClassSigmaP{k-1}$.
Combined with the classical hardness results of Jeroslow and Blair, this establishes $\ClassSigmaP{k-1}$-completeness.
Moreover, we showed that deciding the unboundedness of $k$-level LP is also $\ClassSigmaP{k-1}$-complete, demonstrating that both problems reside at the same level of the polynomial hierarchy. Finally, by extending our analysis to the mixed-binary setting, we proved that the decision version of $k$-level mixed-binary LP is $\ClassSigmaP{k}$-complete, exhibiting a natural one-level increase in complexity when binary variables are introduced.

Taken together, these results yield a precise classification of the complexity of the decision version of $k$-level LP within the polynomial hierarchy.
They highlight the close correspondence between the $k$-level LP and the quantified logical problems.
Future research may investigate whether analogous strong hardness results extend to more general hierarchical optimization models, such as nonlinear or integer formulations, or whether structural or parameterized restrictions on constraints can give rise to tractable subclasses.
Furthermore, it is of interest to investigate the hardness of the decision version of $k$-level LP with any $k$.

% \section*{Acknowledgments}

% This study was funded by the NSERC Grants 2024-04051.

%% file: appendix.tex
\begin{center}
\Large
\underline{\textbf{Online Supplement}}
\end{center}

\section{Challenges to extend Buchheim’s Argument to Higher-Level LPs}
\label{sec:Buchheim}

In recent work~\cite{Buchheim2023}, Buchheim demonstrated the inclusion of the decision version of bilevel LP in \ClassNP{}.
He also briefly discussed that his argument might be extendable to analyze the hardness of the decision version of $k$-level LP for $k \geq 3$ (Remark 2 in the reference).
This section highlights some challenges that arise when attempting such an extension.
We first outline the argument of Buccheim.
Then, we demonstrate the difficulty of the extension of his argument to trilevel LP.

The idea of Buchheim’s argument for bilevel LP can be summarised as follows.
Whenever a bilevel LP instance is feasible, one can identify a dual-feasible basis of the follower’s problem that serves as a polynomial-size certificate of feasibility. 
The validity of this certificate can be checked by solving a corresponding linear system of inequalities, which can be done in polynomial time in the size of the bilevel LP.
This reasoning extends naturally to the decision version of bilevel LP (i.e., the decision of the existence of a feasible solution whose objective value is less than or equal to a given threshold), establishing its inclusion in NP.
A more formal exposition is provided below.

Suppose we have the following bilevel problem, where the leader (the first player) chooses $x_1$ while the follower (the second player) chooses $x_2$:
\begin{align*}
\min_{x_1, x_2} 
\left\{ c_{11}^{\top} x_1 + c_{12}^{\top} x_2 
: 
\begin{array}{l}
A_{11} x_1 + A_{12} x_2 = b_1, x_1 \ge \ZeroVector, \\
x_2 \in \Argmin_{x_2'} \{ c_{21}^{\top} x_1 + c_{22}^{\top} x_2' : A_{21} x_1 + A_{22} x_2' = b_2, x_2' \ge \ZeroVector \} 
\end{array}
\right\}.
\end{align*}
Assume that $A_{22}$ has full row rank and that there exists a feasible solution $(\bar{x}_1, \bar{x}_2)$ whose objective value is less than or equal to $t$, that is,
$$c_{11}^{\top} \bar{x}_1 + c_{12}^{\top} \bar{x}_2 \le t.$$
Consider the lower-level problem with $x_1$ fixed at $\bar{x}_1$.
An optimal basic feasible solution $\hat{x}_2$ can always be found for this lower-level problem, since it has an optimal solution.
Let $B$ denote the index set of basic variables, and let $N = \{1, \ldots, \NVariables_2\} \setminus B$ denote the index set of nonbasic variables.
Reorder the lower-level variables such that $\hat{x}_2 = (\hat{x}_2^B, \hat{x}_2^N)$, where $\hat{x}_2^B$ corresponds to the basic variables.
Partition the constraint matrix and cost vector accordingly as $A_{22} = (A_{22}^B, A_{22}^N)$ and $c_{22} = (c_{22}^B, c_{22}^N)$.
By definition, $A_{22}^B$ is invertible, and the associated basic optimal solution can be written as 
$$
(\hat{x}_2^B, \hat{x}_2^N) = ((A_{22}^B)^{-1} (b_2 - A_{21} \bar{x}_1), \ZeroVector).
$$
The feasibility of $\hat{x}_2$ implies $(A_{22}^B)^{-1} (b_2 - A_{21} \bar{x}_1) \ge \ZeroVector$.
Moreover, since $\bar{x}_2$ is also optimal for the lower-level problem given $x_1 = \bar{x}_1$, it follows that
$$(c_{22}^B)^{\top} (A_{22}^B)^{-1} (b_2 - A_{21} \bar{x}_1) = c_{22}^{\top} \bar{x}_2.$$

A moment of thought gives us that if there is a feasible solution whose objective value is less than or equal to $t$, such a solution can be obtained by nondeterministically selecting a dual-feasible basis $B$, and solving the following system
\begin{align*}
c_{11}^{\top} x_1 + c_{12}^{\top} x_2 &\le t, \\
A_{11} x_1 + A_{12} x_2 &= b_1, x_1 \ge \ZeroVector, \\
A_{21} x_1 + A_{22} x_2 &= b_2, x_2 \ge \ZeroVector, \\
(c_{22}^B)^{\top} (A_{22}^B)^{-1} (b_2 - A_{21} x_1) &= c_{22}^{\top} x_2, \\
(A_{22}^B)^{-1} (b_2 - A_{21} x_1) &\ge \ZeroVector.
\end{align*}
A feasible solution with objective value at most $t$ exists if and only if there is a nondeterministic choice of dual-feasible basis $B$ for which the above system is consistent.

We now turn our attention to multilevel LP.
Buchheim discussed that this reasoning might extend to establish that the decision version of the $k$-level LP lies in the class~\ClassSigmaP{k-1}.
He suggested applying an analogous reduction, in which the decision version of a $k$-level LP is reformulated as a feasibility question for a $(k-1)$-level LP.
He briefly outlined his argument to replace the optimality condition of the final follower with a system of linear inequalities.
While this idea is conceptually appealing, we believe that additional assumptions would be required for such an extension to be valid.
The following counterexample shows that a naive application of this reasoning can lead to an incorrect conclusion. 
To be more specific, nondeterministically guessing a dual-feasible basis of the last follower’s problem and following the above reasoning for bilevel LP may falsely indicate that a trilevel LP instance is feasible.

Consider the problem
\begin{equation}
\min_{x_1}
\max_{x_2}
\min_{x_3}
\{ x_3 : x_3 \ge x_2 - 1, x \ge \ZeroVector \}.
\label{eq:Buchheim}
\end{equation}
The third player selects $x_3 = \max\{x_2 - 1, 0\}$, so that the second player’s problem becomes $\max_{x_2} \{ \max\{x_2 - 1, 0\} : x_2 \ge 0 \}$, which is unbounded.
Therefore, problem~\eqref{eq:Buchheim} is infeasible.

Suppose we are interested in the existence of a feasible solution with objective value at most $t$ for some nonnegative rational number $t$.
Clearly, no such solution exists, since~\eqref{eq:Buchheim} is infeasible.
Following Buchheim’s line of reasoning, we may express the third player’s problem in standard form as
$$
\min_{x_1}
\max_{x_2}
\min_{x_3, s_3}
\{ x_3 : x_3 - s_3 = x_2 - 1, x \ge \ZeroVector, s_3 \ge 0 \}.
$$
Now select the basic solution to the third player’s problem corresponding to the basic variable $s_3$, which is dual feasible. 
Applying Buchheim’s approach, we replace the third player’s problem by the constraints enforcing feasibility of this basic solution:
$$
\min_{x_1}
\max_{x_2, x_3, s_3}
\{ x_3 : x_3 \le t, x_3 - s_3 = x_2 - 1, x \ge \ZeroVector, s_3 \ge 0, x_3 = 0, 1 - x_2 \ge 0 \}.
$$
This bilevel LP is feasible, as $x = \ZeroVector$ and $s_3 = 1$ satisfy all constraints as well as the optimality of the second player's problem.
Hence, if one follows the argument literally, one would erroneously conclude that~\eqref{eq:Buchheim} admits a feasible solution with objective value at most $t$.
Although the idea of reducing decision problems involving $k$-level LPs to those of $(k-1)$-level LPs using simple linear inequalities is theoretically appealing, it appears to require additional assumptions to ensure its correctness.
A thorough investigation of the conditions under which such an extension becomes valid is, however, beyond the scope of this work.